\documentclass{article}
\usepackage{psfrag,graphicx}

\newtheorem{thm}{Theorem}[section]

\newtheorem{cor}[thm]{Corollary}
\newtheorem{prop}[thm]{Proposition}
\newtheorem{defn}[thm]{Definition}
\newtheorem{example}[thm]{Example}
\newtheorem{rmk}[thm]{Remark}
\newtheorem{question}[thm]{Question}

\newcommand{\bC}{{\rm\bf C}}
\newcommand{\bR}{{\rm\bf R}}
\newcommand{\bZ}{{\rm\bf Z}}
\newcommand{\cC}{{\cal C}}
\newcommand{\cL}{{\cal L}}

\newcommand{\bdy}{\partial}
\newcommand{\del}{\partial}
\newcommand{\wtilde}{\widetilde}
\newcommand{\what}{\widehat}

\newcommand{\eps}{\epsilon}
\newcommand{\Sig}{{\Sigma}}
\newcommand{\om}{\omega}
\newcommand{\vphi}{\varphi}

\newcommand{\ol}{\overline}

\def\sqr#1#2{{\vcenter{\hrule height.#2pt
   		\hbox{\vrule width.#2pt height#1pt \kern#1pt
      		\vrule width.#2pt}\hrule height.#2pt}}}
\def\square{\mathchoice\sqr67\sqr67\sqr{2.1}6\sqr{1.5}6}

\newcommand{\proof}[1]{\noindent{\bf Proof#1:\  }}

\title{Symplectic rational blowdowns}

\author{Margaret Symington 
\thanks{The author is grateful for the support of an NSF post-doctoral 
fellowship, DSM9627749.}
\\
{\small Department of Mathematics}\\  
{\small University of Texas, Austin, TX 78712}\\
{\small msyming@math.utexas.edu} }

\begin{document}

\maketitle

\begin{abstract}
We prove that the rational blowdown, a surgery on smooth $4$-manifolds
introduced by Fintushel
and Stern, can be performed in the symplectic category.  As a consequence, 
interesting families of smooth $4$-manifolds, including the exotic $K3$ 
surfaces
of Gompf and Mrowka, admit symplectic structures.
\end{abstract}

A basic problem in symplectic topology is to understand what
smooth manifolds admit a symplectic structure (a closed
non-degenerate $2$-form).  In this paper we focus on this question in
dimension $4$.
Currently, the primary methods for constructing
smooth (irreducible) $4$-manifolds in such a way that one can distinguish
them by Donaldson or Seiberg-Witten invariants
are surgery constructions that use complex manifolds as building blocks.  
These surgery methods are (smooth) logarithmic transforms, 
rational blowdowns, and connect sums along surfaces.
It is interesting to see when these surgeries can be performed in the 
symplectic category.
In this paper we prove that performing a rational blowdown of a symplectic
manifold along symplectic surfaces yields a symplectic manifold.
This result establishes that certain exotic $4$-manifolds,
including the exotic $K3$ surfaces of Gompf and Mrowka~\cite{GM}, 
are symplectic.

In any even dimension, two symplectic manifolds can be summed
along codimension $2$ symplectic submanifolds to yield a symplectic 
manifold.  
We refer to this symplectic operation, which was proposed by Gromov~\cite{Gr2},
as the symplectic sum.  
Gompf~\cite{Go} used the symplectic sum to construct a plethora of interesting
symplectic manifolds, including
the first examples of simply connected symplectic $4$-manifolds 
that are not homotopic to any complex surface and some exotic $K3$ surfaces.
More recently, Fintushel and Stern~\cite{FS2} have used the connect
sum along smoothly embedded tori to produce a rich class of
exotic $4$-manifolds, some homeomorphic to a $K3$ surface, 
many of which cannot admit a symplectic structure.

\begin{figure} 
\begin{center}
        \psfragscanon
        \psfrag{np2}{$n+2$}
        \psfrag{2}{$2$}
        \includegraphics[width=2in]{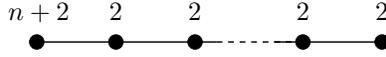}
\end{center}
\caption{The plumbing diagram for $C_n$, $n\ge2$}
\label{Cn.fig}
\end{figure}

The logarithmic transform was first studied in the smooth category by
Gompf and Mrowka~\cite{GM} who used it to produce the first examples
of irreducible $4$-manifolds that are not complex.
Subsequently, Fintushel and Stern~\cite{FS} introduced the rational
blowdown and showed that in certain situations a smooth
logarithmic transform can be achieved via a sequence of blowups followed by 
a rational blowdown. 
The rational blowdown is a surgery in which a neighborhood of a chain of 
spheres $C_n$, $n\ge2$, represented by the plumbing diagram in
Figure~\ref{Cn.fig}
is replaced by a rational (homology) ball $B_n$. 

Using the rational blowdown, Fintushel and Stern 
constructed other interesting examples of smooth $4$-manifolds.
Their examples led them to ask whether the rational blowdown of a
symplectic $4$-manifold along symplectic spheres is a symplectic operation.
Theorem~\ref{blowdown.thm} of the next section asserts the answer is yes.
As a consequence, an infinite family of surfaces not homotopic
to a complex surface, constructed by Fintushel and Stern in~\cite{FS}, are
symplectic.  Furthermore, the complete 
set of Gompf-Mrowka examples of exotic $K3$ surfaces~\cite{GM} 
are also symplectic, extending the results in~\cite{Go}.

The Fintushel-Stern examples that are symplectic as a consequence of
Theorem~\ref{blowdown.thm} are constructed from 
the simply connected minimal elliptic surfaces $E(n)$, $4\le n\in\bZ$,
that have Euler characteristic $\chi(E(n))=12n$.
In $E(n)$ one can find two copies of $C_{n-2}$, embedded so that
the spheres are symplectic (consult~\cite{Fu} and~\cite{FS}).
Performing a symplectic rational blowdown along one of these chains of
spheres yields a manifold $G(n)$ whose homotopy type is different
from any complex manifold.
Blowing down both chains of spheres yields a manifold diffeomorphic to
a Horikawa surface $H(n)$, a complex manifold of general type 
(which is therefore K\"ahler).

The Gompf-Mrowka examples
are obtained from the $K3$ surface by performing
smooth logarithmic transforms along three pairs of tori
(in which a neighborhood of a torus having trivial normal bundle
is removed and replaced using
a diffeomorphism of the boundary not homotopic to the identity).
The resulting manifolds are denoted $K(p_1,q_1;p_2,q_2;p_3,q_3)$, with
$p_i,q_i$ relatively prime.
Modulo certain relations between the $p_i,q_i$,
these manifolds are mutually non-diffeomorphic but are
homeomorphic to either a $K3$ surface or $3\bC P^2\# 19\overline{\bC P}^2$
In~\cite{Go}, Gompf showed that $K(p_1,q_1;1,1;p_3,q_3)$ are symplectic
by presenting them as symplectic sums of simply connected Dolgachev surfaces.
The work of Fintushel and Stern~\cite{FS} shows that all of the
$K(p_1,q_1;p_2,q_2;p_3,q_3)$ can be constructed by blowing up points
and then performing a rational blowdown.
Because the necessary submanifolds can be chosen to be symplectic, 
Theorem~\ref{blowdown.thm} implies that all the $K(p_1,q_1;p_2,q_2;p_3,q_3)$ 
are symplectic.

In the next section we give a precise definition of the
symplectic rational blowdown and state the main theorem.
The essence of the proof of Theorem~\ref{blowdown.thm} is in
our choice of model spaces for $C_n$ and a collar neighborhood of
the boundary of $B_n$.  
Indeed, using a model for $L(n^2,n-1)\times (0,\infty)$ as a guide, we endow
$B_n$ with a symplectic structure such that the complement of the spheres
in $C_n$ is symplectomorphic to a collar of $B_n$.
The gluing is then clear.
To describe the model spaces we use symplectic boundary reduction which
is the main step in the 
procedure of symplectic cutting (defined by Lerman~\cite{Ler}).
We define symplectic boundary reduction in Section~\ref{bdy.sec} and
construct our model spaces in Section~\ref{model.sec}.
We then prove Theorem~\ref{blowdown.thm} in the last section.

\begin{rmk} \label{proofs.rmk} \rm
Theorem~\ref{blowdown.thm} can also be deduced as a straightforward
application of the $3$-fold sum, a 
symplectic surgery developed by the author.  
The $3$-fold sum is a sum along positively intersecting symplectic 
surfaces that is part of 
a generalization (in dimension $4$) of the symplectic sum~\cite{Sy1}.
We sketch this alternative proof of Theorem~\ref{blowdown.thm}
in Remark~\ref{3sum.rmk}, referring the reader to~\cite{Sy2} for
details on the $3$-fold sum.  
\end{rmk}

We use the notation $[b_1,b_2,\ldots, b_n]$ to denote the {\it negative}
continued fraction expansion 
$b_1 - 1/(b_2 - 1/(\ldots - 1/ b_n )\ldots)$.

\section{The symplectic rational blowdown} \label{blowdown.sec}

The symplectic
rational blowdown  generalizes the blowing down of a $-4$ sphere (a
sphere with self-intersection $-4$), in
which a neighborhood of the sphere is replaced by the complement
of a conic in $\bC P^2$.  As observed by Gompf~\cite{Go}, this can be 
achieved using the
symplectic sum (assuming the $-4$ sphere is symplectic).

Let $C_n$, $n\ge2$, be a tubular neighborhood of a union of spheres 
$\cup_{i=1}^{n-1} S_i$
such that $S_1\cdot S_1= -(n+2)$, $S_i\cdot S_i=-2$ for $i=2,\ldots n-1$,
$S_i\cdot S_{i+1} = 1$, for $i=1,\ldots n-2$, and $S_i\cdot S_j=0$ otherwise.
Thus $C_n$ is a plumbing of disk bundles over spheres represented by
the diagram in Figure~\ref{Cn.fig}.
The boundary of $C_n$ is the lens space $L(n^2,n-1)$ which also bounds a
manifold $B_n$ that has the same rational homology as a ball (see~\cite{CH}).
In Section~\ref{model.sec} we define symplectic models $(C_n,\om_{C_n})$
and $(B_n,\om_{B_n})$ whose symplectic structures depend on the areas of the
spheres $\{S_i\}_{i=1}^{n-1}$.  
The symplectic structure $\om_{C_n}$ is chosen so that the spheres
$\{S_i\}_{i=1}^{n-1}\subset(C_n,\om_{C_n})$ are symplectic and intersect
orthogonally with respect to $\om_{C_n}$.
(Note that an
orthogonal intersection of symplectic surfaces is necessarily positive.)

\begin{defn}
Suppose there is a symplectic embedding 
$\psi:(C_n,\om_{C_n})\rightarrow (M,\om)$.
Let $M^- = M-\psi\left(\cup_{i=1}^{n-1}S_i\right)$ and
let $B_n$ be a rational homology ball with no prescribed symplectic
structure.
A {\bf symplectic rational blowdown} of $(M,\om)$ along the spheres
$\psi\left(\cup_{i=1}^{n-1}S_i\right)$ is a closed manifold 
$\wtilde M = M^-\cup_{\vphi} B_n$
with a symplectic structure $\tilde\om$ such that
$(M^-,\tilde\om)$ and $(M^-,\om)$ are symplectomorphic.
\end{defn}
 
\begin{rmk} \label{positive.rmk} \rm
If $(M,\om)$ contains a union of symplectic spheres with the same intersection
pattern as $\cup_{i=1}^{n-1}S_i\subset C_n$, then there is a symplectic embedding
$\psi:(C_n,\om_{C_n})\hookrightarrow (M,\om)$.  
In fact, the spheres in $M$ can be isotoped 
to make the intersections orthogonal, keeping them symplectic all the 
while (cf.~\cite{McDP}).  
Then by a version of the symplectic neighborhood theorem 
(Proposition~\ref{nbhd.prop})  $C_n$ is symplectomorphic to a neighborhood
of the isotoped spheres.
\end{rmk}

\begin{thm} \label{blowdown.thm}
Suppose there is a symplectic embedding
$\psi:(C_n,\om_{C_n}) \hookrightarrow (M,\om)$.
Then there exists a symplectic rational ball $(B_n,\om_{B_n})$ 
and a symplectic map $\vphi$ making
$\wtilde M = M^- \cup_{\vphi} (B_n,\om_{B_n})$
a symplectic rational blowdown of $M$.
The volume of $\wtilde M$ is determined by the volume of $M$ and
the areas of the spheres $\{S_i\}_{i=1}^{n-1}$.
\end{thm}
Note that a smooth rational blowdown is well defined up to diffeomorphism  
because any diffeomorphism of
the boundary of $B_n$ extends over the rational ball~\cite{Bon}.
The symplectic blowdown of a $-1$ sphere is unique because
the symplectic structure of any ball that is standard near the boundary
is diffeomorphic to the standard structure via a diffeomorphism
that is the identity near the boundary~\cite{Gr1}.  
It is an interesting question
whether a symplectic rational blowdown is also unique up to 
symplectomorphism.

\section{Symplectic boundary reduction} \label{bdy.sec}

Let $(M,\om)$ be a symplectic $4$-manifold whose boundary is a circle bundle
over a surface $\Sig$.  
Suppose that all vectors tangent to the circle fibers lie in the kernel 
of $\om|_{\bdy M}$.
Then there is a closed symplectic manifold $(\what M,\hat\om)$, 
unique up to symplectomorphism, that contains an embedded copy of $\Sig$
and is such that $(M-\bdy M,\om)$ and $(\what M-\Sig,\hat\om)$ are
symplectomorphic.
We call $(\what M,\hat\om)$ the {\bf symplectic boundary reduction} of 
$(M,\om)$.
It can be realized as the image of a map $\pi$ which is symplectic when
restricted to the interior of $M$,
and which collapses each circle fiber of $\bdy M$ to a point.
The image $\pi(\bdy M)$ is the embedded copy of $\Sig$;
it is a symplectic submanifold of $\what M$.
(The above description is true in higher dimensions with $\Sig$ being
a manifold of dimension $2$ less than that of $M$).

In a neighborhood of a fiber
of $\bdy M$, the map $\pi$ can be described in local coordinates 
as follows.
Any fiber in the boundary of $M$ has a neighborhood symplectomorphic to
$(D^2\times A^2,dx_1\wedge dy_1 + dp_2\wedge dq_2)$ where 
$A^2=\{0\le p_2 <\eps\} \subset \bR\times S^1$ and $q_2$ is defined mod 1.
With respect to these local coordinates, 
$\pi$ is the projection 
\[\pi:\ 
(x_1,y_1,p_2,q_2)\to \left(x_1,y_1,\sqrt{\frac{p_2}{\pi}}\cos(2\pi q_2),
\sqrt{\frac{p_2}{\pi}}\sin(2\pi q_2)\right).
\]

One can also take the symplectic boundary reduction of a symplectic manifold
when its boundary is not smooth, but rather has corners.
Specifically, we allow the boundary of $(M,\om)$ to
have more than one smooth component, pairs of which meet along
Lagrangian tori (tori $T$ of half the dimension of $M$ such that
$\om|_T=0$).  
The definition of symplectic boundary reduction in this context
is the same as above except
that the interior of $M$ is symplectomorphic to the complement of
a union of intersecting symplectic surfaces in $\what M$.
Examples~\ref{coord.ex} and~\ref{cone.ex} are local models for 
boundary reduction near a corner on the boundary of $M$.  
Note that we always take the boundary reduction only along
the closed part of $\bdy M$.

Here and throughout this paper we use models that are obtained
from $T^*T^2=\bR^2\times T^2$ with the standard symplectic structure
$\om_0=dp\wedge dq$ where $p=(p_1,p_2)$ are coordinates on $\bR^2$
and $q=(q_1,q_2)$ are coordinates on $T^2$ defined mod 1.

\begin{example}  \label{coord.ex} \rm
{\bf : $\mathbf{(\bR^4,dx\wedge dy).}\quad$}
Let $Q$ be the first quadrant of $\bR^2$ and $\ol Q$ its closure.
Consider $Q\times T^2 \subset (T^*T^2,\om_0)$ and
define the map $\pi:Q\times T^2\rightarrow \bR^4$ with coordinates
$(x_1,y_1,x_2,y_2)$ by the formula
\[
(x_i,y_i)=\left(\sqrt{\frac{p_i}{\pi}}\cos(2\pi q_i),
\sqrt{\frac{p_i}{\pi}}\sin(2\pi q_i)\right).
\]
It is a symplectomorphism between $Q\times T^2$ and the complement of 
the coordinate planes $x_1=y_1=0$ and $x_2=y_2=0$ in
$(\bR^4,dx\wedge dy)$.
Extending $\pi$ to $\ol Q\times T^2$ we get a projection
to $\bR^4$ in which the image of the torus $p_1=p_2=0$ is the origin
and the image of each circle fiber on the rest of the boundary of
$\ol Q\times T^2$ is a point on one of the coordinate planes.
The image of this projection, which is all of $\bR^4$,
is the boundary reduction of $\ol Q\times T^2$.
\end{example}

\begin{example} \label{cone.ex} \rm
{\bf : $\mathbf{(\bR^4,dx\wedge dy)\ }$ again.\ \ }
Now consider any closed positive cone $C$ in $\bR^2$ defined by integral 
vectors $u$ and $v$ such that the matrix $B=[u\, v]$ is in $GL(2,\bZ)$.
The boundary reduction of 
$C\times T^2\subset (T^*T^2,\om_0)$ is $\bR^4$ with the standard
symplectic structure.
This follows because $(C\times T^2,\om_0)$ is symplectomorphic to 
$(\ol Q\times T^2,\om_0)$ via the
map $\vphi(p,q)=(Bp+r,B^{-T}q)$ where $r$ is the vertex of the cone 
$C\subset \bR^2$.
\end{example}

\begin{defn} \label{along.defn}
If $\Sig$ is a surface in the image of $\bdy M$ under 
symplectic boundary reduction of
$(M,\om)$, then 
we call the preimage $\pi^{-1}(\Sig)$ the {\bf boundary along $\Sig$}.
\end{defn}
 
\section{Model neighborhoods} \label{model.sec}

By construction, all of our model spaces admit Hamiltonian $2$-torus
actions, and all of the figures we draw in $\bR^2$ are in fact images
of moment maps for the torus actions, though we do not appeal to that 
language in this paper.  However, the reader should note that this means that
in our figures, each interior point represents a torus, each edge point
represents a circle and each vertex represents a single point.

We begin by describing a symplectic structure $\om_{n,m}$
on $V_{n,m}=L(n,m)\times(0,\infty)$, for $n\ge m\ge 1$ relatively prime,
that is induced
from the standard structure $\om_0=dp\wedge dq$ on $T^*T^2$
via boundary reduction.
Recall that a lens space $L(n,m)$ can be presented as 
the union of two solid tori glued together via a map $\phi$ 
of their boundaries such that $\phi_*\mu_2= -m\mu_1 + n\lambda_1$ 
where $\mu_i,\lambda_i$ are meridinal and longitudinal cycles 
on the boundaries.

\begin{example} \label{Lnm.ex} \rm
{\bf : }{\boldmath$(V_{n,m},\om_{n,m}).\quad$}
Let $(V_{n,m},\om_{n,m})$ be the boundary reduction
of $U_{n,m}\times T^2 \subset (T^*T^2,\om_0)$ where 
$U_{n,m}=\{p_1\ge 0\}\cap\{p_2\ge \frac{m}{n}p_1\}\cap\{p_2>0\}$.
The $3$-dimensional submanifold of $V_{n,m}$ 
that is the image of $(\{p_2=1\}\cap U_{n,m})\times T^2$ is a union of 
two solid tori with shared boundary
$\{p_1=c,\ p_2=1\}\times T^2$ for some $0<c<\frac{m}{n}$.  
The boundary of the solid torus that is the image of 
$(\{p_1\le c,\ p_2=1\}\cap U_{n,m})\times T^2$ 
has a meridian whose tangent vectors are
$\frac{\del}{\del q_1}$, while the other one has a meridian whose tangent
vectors are $-m\frac{\del}{\del q_1}+n\frac{\del}{\del q_2}$.
Since these two meridians are identified, the $3$-manifold is 
$L(n,m)$.  
Hence, $(V_{n,m},\om_{n,m})$ is a symplectic model for
$L(n,m)\times(0,\infty)$.
\end{example}

For the following examples, let $L_0,L_1$ be the lines $\{p_1=0\}$,
$\{p_2=0\}$ in $\bR^2$.

\begin{example} \label{one.ex} \rm
{\bf :}{\boldmath $(N_b,\om_{N_b}).\quad$}
Let $L_2$ be the line $\{p_2=\frac{1}{b}(p_1-a)\}$, $a>0$.
Consider the closed domain in $\bR^2$ lying between $L_0$ and $L_2$
and bounded below by $L_1$.
Let $U_b$ be a neighborhood of $L_1$ in this domain.
Then the boundary reduction $(N_b,\om_{N_b})$ 
of $U_b\times T^2$ (taken along the
closed edges) is a symplectic neighborhood of a sphere of 
self-intersection $-b$ and area $a$.  
Indeed, the image of $(L_1\cap U_b)\times T^2$ 
is a sphere of area $a$ since it is
the boundary reduction of a cylinder
$\{(p_1,q_1)|0\le p_1\le a\}$ with symplectic form $dp_1\wedge dq_1$.
As per Examples~\ref{coord.ex} and~\ref{cone.ex}, $N_b$ 
is the union of two polydisks, and the boundary of $N_b$ is $L(b,1)$,
verifying that the sphere has self-intersection $-b$.
\end{example}

\begin{example}  \label{nbhd.ex} \rm
{\bf :}
{\boldmath $(C_n,\om_{C_n})$ and $(C^-_n,\om_{C^-_n}).\quad$}
Suppose the spheres $\{S_i\}_{i=1}^{n-1}$
have areas $\{a_i\}_{i=1}^{n-1}$.
Let $n_i,m_i$ be the relatively prime integers such that 
$\frac{n_i}{m_i}=[n+2,2,\ldots 2]$ for a continued fraction expression
of length $i$, and define
vectors $r_i = \left(\begin{array}{c}n_{i-1} \\ m_{i-1}\end{array}\right)$ 
with $n_0=1$ and $m_0=0$.

For $i=2,\ldots n$, let $t_i=\sum_{j=1}^{i-1}a_jr_j$ and define
lines $L_i$, $i=2,\ldots,n$ by the parametric equations 
$ t r_i + t_i$, $t\in(-\infty,\infty)$.
Consider the closed domain bounded below by the union of lines
$L_i$, $i=1,\ldots, n-1$ and lying between the lines $L_0,L_{n}$.
Let $U_{C_n}$ be a neighborhood in this domain of the lines 
$L_i$, $i=1,\ldots, n-1$.
The boundary reduction of $U_{C_n}\times T^2$ is a symplectic model
for $(C_n,\om_{C_n})$.  

To see this notice that $U_{C_n}\times T^2$ is the union of 
$U_{n+2}\times T^2$ (as defined in Example~\ref{one.ex}) 
and the images of $n-2$ copies of $U_{2}\times T^2$ under symplectic maps 
$(p,q)\mapsto (T_ip+t_i,T_i^{-T}q)$, $i=2,\ldots n-1$ where
$T_i=R_{n+2}R_2^{i-2}$ with
\[
R_{k}=\left( \begin{array}{cc} k & -1 \\ 1 & 0 \end{array} \right).
\]
These maps induce symplectomorphisms that plumb the disk bundles 
$N_{n+2},$ $N_2,\ldots N_2$ to form $(C_n,\om_{C_n})$.  

We let $(C^-_n,\om_{C^-_n})$ be the boundary
reduction of $U_{C^-_n}\times T^2$ where
$U_{C^-_n}=U_{C_n}-\{L_i\}_{i=1}^{n-1}$.
This is symplectomorphic to the complement of the surfaces 
$\{S_i\}_{i=1}^{n-1}\subset (C_n,\om_n)$ and hence
is a symplectic model for a collar neighborhood of the boundary of $M^-$.
\end{example}

We now define a symplectic rational ball $(B'_n,\om_{B'_n})$ that
is a complement of two spheres in a rational ruled surface.
In the proof of Theorem~\ref{blowdown.thm} we will see how to modify the
symplectic structure on
this rational ball to obtain the ball $(B_n,\om_{B_n})$ which is required
for the gluing.
Let $F_{n-1}$ be a rational ruled surface 
that contains symplectic sections
$\Sig_{n+1}$ and $\Sig_{-n+1}$ with self-intersections $n+1$,$-n+1$.  
(For instance $F_{n-1}$ can be a projectivized plane bundle
with holomorphic sections $\Sig_+$ and $\Sig_-$ of 
self-intersections $n-1$,$-n+1$ respectively.  
In this case $[\Sig_{n+1}]=[\Sig_+]+[f]$ and $[\Sig_{-n+1}]=[\Sig_-]$
where $f$ is  a fiber.)
Since $\Sig_{n+1}\cdot\Sig_{-n+1}=1$ and the spheres span 
the rational homology of $F_{n-1}$, 
the complement $F_{n-1}-(\Sig_{n+1}\cup\Sig_{-n+1})$ is a rational ball
with boundary $L(n^2,n-1)$.
Let $(B'_n,\om_{B'_n})$ be this rational ball with the symplectic structure 
inherited from $F_{n-1}$.
(Note that the ruled surface $F_{n-1}$ has a symplectic structure, 
well-defined up to symplectomorphism
by the areas $\alpha_{n+1},\alpha_{-n+1}$
of the two sections~\cite{McD4}.)

\begin{example} \label{ball.ex} \rm
{\bf :} {\boldmath $(A'_n,\om_{A'_n}).\quad$}
We can assume that the two sections $\Sig_{n+1},\Sig_{-n+1}$ intersect
orthogonally with respect to the symplectic structure on $F_{n-1}$,
isotoping one of them if necessary (cf.~\cite{McDP}).
Define $L_2$, $L_3$ 
parametrically in $t$ by $ t r_2 + \alpha_{n+1} r_1$ and 
$t r_3 + \alpha_{-n+1} r_2 + \alpha_{n+1} r_1$ where
\[
r_1 = \left(\begin{array}{c}1\\0\end{array}\right), \qquad
r_2 = \left(\begin{array}{c} -n-1 \\ -1 \end{array}\right) \quad {\rm and} 
\quad
r_3 = \left(\begin{array}{c} -n^2 \\ -n+1 \end{array}\right).
\]
Assuming $\alpha_{n+1}>(n+1)\alpha_{-n+1}$,
consider the closed domain in $\bR^2$ that lies between $L_0$, $L_3$ and
below $L_1$, $L_2$.
Let $U_{A'_n}$ be a neighborhood of $L_1\cup L_2$ in this domain, minus
the two lines.
The boundary reduction $(A'_n,\om_{A'_n})$ of $U_{A'_n}\times T^2$ is a 
model for a collar neighborhood of the boundary of $(B'_n,\om_{B'_n})$.
\end{example}

For our model neighborhoods to be useful 
we need the following simple modification of the
symplectic neighborhood theorem.

\begin{prop} \label{nbhd.prop}
Consider two embeddings of a configuration of spheres,
$j:\cC\hookrightarrow(M,\om)$ and
$j':\cC\hookrightarrow(M',\om')$, such the areas and
self-intersections of the images of a given sphere are
the same and all intersections are orthogonal (and hence positive)
with respect to the ambient symplectic form.  
Then $j(\cC)$ and $j'(\cC)$ have symplectomorphic neighborhoods.
\end{prop}

The proof is a standard application of Moser's
method to turn a diffeomorphism into a symplectomorphism except that
one must be careful in the vicinity of the intersection points.
Details of how to handle such situations were provided in~\cite{McR}.

\begin{cor} \label{bdy.cor}
Suppose $(M,\om)$ and $(M',\om')$ are symplectic manifolds with boundary
such that the kernel of the symplectic form (restricted to the boundary)
defines a union of smooth components fibered by circles and intersecting
along Lagrangian tori.  
Let $\pi(M),\pi'(M')$ be their symplectic boundary
reductions.  If $\pi(\bdy M),\pi'(\bdy M')$ define configurations of
symplectomorphic submanifolds with the same intersection patterns in
$\pi(M),\pi'(M')$, then $\bdy M, \bdy M'$ have symplectomorphic
collar neighborhoods.
\end{cor}

\proof{}
This follows from the fact that there is a unique way to extend
to the boundaries of $M,M'$ the symplectic map
$(\pi')^{-1}\circ\phi\circ\pi$ where $\phi$ is a symplectomorphism
of neighborhoods of $\pi(\bdy M)$ and $\pi'(\bdy M')$.
\hfill$\square$

\section{Proof of Theorem~\ref{blowdown.thm}}

\proof{}
To blow down a symplectic $-4$-sphere in $(M,\om)$ one simply takes
the symplectic sum of $M$ with $\bC P^2$ along the $-4$-sphere in $M$ and
a conic $Q$ in $\bC P^2$, as shown in~\cite{Go}.
This is the rational blowdown for $n=2$.
For notational convenience we restrict our attention to the cases $n\ge 3$.

Recall the domains $U_{n,m}$, $U_{C^-_n}$ and $U_{A'_n}$ defined in
Examples~\ref{Lnm.ex}, \ref{nbhd.ex}, and~\ref{ball.ex}.  
The boundary reduction of the product of each of these with
$T^2$, viewed as a subset of $(T^*T^2,\om_0)$, is a symplectic model
for $V_{n,m}$, $C^-_n$ and $A'_n$ respectively.
Corollary~\ref{bdy.cor}
implies that there are symplectic embeddings 
$\psi_1:(C^-_n,\om_{C^-_n})\hookrightarrow M-\psi(\cup_{i=1}^{n-1}S_i)$ and
$\psi_2:(A'_n,\om_{A'_n})\hookrightarrow (B'_n,\om_{B'_n}) 
= F_{n-1}-\{\Sig_{-n+1},\Sig_{n+1}\}$.

Notice that there is a translation of
the domain $U_{C^-_n}$ that is a subset of $U_{n^2,n-1}$ such that
its closed edges are subsets of the two edges of $U_{n^2,n-1}$.
Since translation in the $p$-coordinates is a symplectomorphism of
$T^*T^2$, this implies that $(C^-_n,\om_{C^-_n})$ is symplectomorphic
to a submanifold of $(V_{n^2,n-1},\om_{n^2,n-1})$.
Call this symplectic embedding $\phi_1$.

Now choose a rational surface $F_{n-1}$ that has sections 
$\Sig_{-n+1},\Sig_{n+1}$
with areas $\alpha_{n+1}>(n+1)\alpha_{-n+1}>0$ such that 
\[
(n-1)\alpha_{n+1}+\alpha_{-n+1} <  
      \sum_{i=1}^{n-1} \left( (n-1)n_{i-1} - n^2 m_{i-1}\right)a_i
\]
where the $n_i,m_i$ are defined in Example~\ref{nbhd.ex} and
the $a_i$ are the areas of the spheres $S_i\subset C_n$.
Because $\alpha_{n+1}>(n+1)\alpha_{-n+1}>0$ there is a
symplectic embedding $\phi_2:A'_n\hookrightarrow V_{n,m}$;
like $\phi_1$ it is a translation in the $p$-coordinates.
Because of our choice of areas, the image of $U_{A'_n}$ under $\phi_2$
lies below the image of $U_{C^-_n}$ under $\phi_1$.
Hence, the union $\phi_1(C^-_n)\cup\phi_2(A'_n)$ is a collar neighborhood of
a submanifold $A_n$ of $V_{n^2,n-1}$ that can be simultaneously glued
onto $(M-\psi(\cup_{i=1}^{n-1}S_i),\om)$ and $(B'_n,\om_{B'_n})$
via symplectomorphisms $\psi_1\circ\phi_1^{-1}$ and $\psi_2\circ\phi_2^{-1}$.
Letting $B_n = 
B'_n \cup_{\psi_2\circ\phi_2^{-1}} A_n$ 
with induced symplectic structure  $\om_{B_n}$, the manifold
\[
\wtilde M = 
\left(M-\psi(\cup_{i=1}^{n-1}S_i)\right)\cup_{\psi_1\circ\phi_1^{-1}} B_n
\]
with the induced symplectic structure $\tilde \om$ is a symplectic rational 
blowdown of $(M,\om)$ along the spheres $\psi(\cup_{i=1}^{n-1}S_i)$. 

\begin{figure} 
\begin{center}
        \psfragscanon
        \psfrag{An}[][]{$A'_4$}
        \psfrag{Cnm}[][]{$C_4^-$}
        \psfrag{Vnm}[][]{$V_{16,3}$}
	\includegraphics[width=3in]{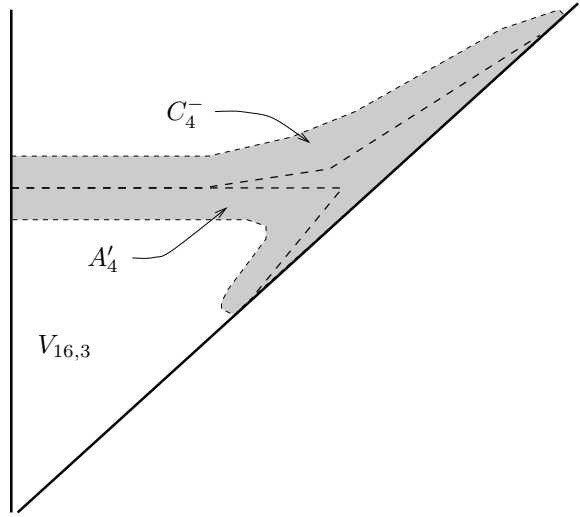}
\end{center}
\caption{Images of $C_4^-$ and $A'_4$ in $V_{16,3}=L(16,3)\times (0,\infty)$}
\label{fit.fig}
\end{figure}

Figure~\ref{fit.fig} illustrates symplectic 
embeddings $\phi_1,\phi_2$ 
of $C_4^-$ and $A'_4$ into $V_{16,3}=L(16,3)\times(0,\infty)$
by showing the images of $U_{C^-_4}$ and $U_{A'_4}$ in $U_{n,m}$. 
Here we have chosen 
$3\alpha_5 + \alpha_{-3}=\sum_{i=1}^3(3n_{i}-16m_{i})a_i$.
The shaded domain represents $A_4$, 
a collar neighborhood of the boundary of the ball $B_4$ we need for  
gluing.

\begin{figure} 
\begin{center}
        \psfragscanon
        \psfrag{Cnm}{$C_4^-$}
        \psfrag{Z}{$Z$}
        \includegraphics[width=3in]{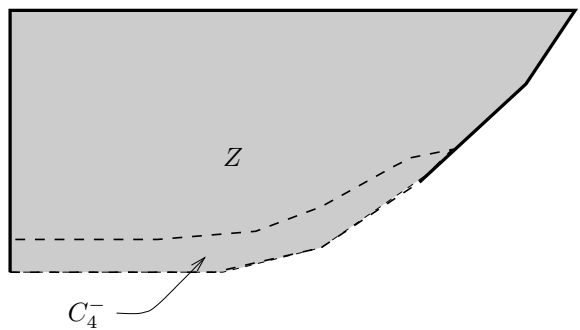}
\end{center}
\caption{Completion of rational balls $B_{4,i}$, $i=1,2$ to $F_{3}$}
\label{ratl.fig}
\end{figure}

The volume of $\wtilde M$ is independent of any choice of rational
ball that fits.
Indeed, suppose $B_{n,1},B_{n,2}$ are two choices of rational balls
for which a symplectic rational blowdown exists.  
By shrinking the gluing loci, we can assume that the boundaries of
the $B_{n,i}$ both have collar neighborhoods that are symplectomorphic
to a common $C_n^-$.
But $C_n^-$ is also a collar neighborhood of the 
boundary of a neighborhood $Z$ of spheres 
$\Sig'_{n+1}\cup\Sig'_{-n+1}$
of self intersections $n+1, -n+1$.  (See Figure~\ref{ratl.fig}.)
Therefore for some symplectic map $\vphi$, each $B_{n,i}\cup_\vphi Z$ 
is a closed symplectic manifold containing a sphere of positive 
self-intersection, and hence is a rational ruled surface~\cite{Gr1}.
The cohomology class of the symplectic form on each of these ruled
surfaces is set by the areas of the spheres $\Sig'_{n+1},\Sig'_{-n+1}$,
thus the volume is the same in both cases.
But this implies that the volume of $B_{n,i}$ is independent of $i$.
\hfill$\square$

Note that in contrast to case of blowing down a $-4$-sphere, 
collar neighborhoods $(C_n^-,\om_{C_n^-})$ and $(A'_n,\om_{A'_n})$ 
cannot be symplectomorphic for
$n\ge 3$ so long as we choose $(B'_n,\om_{B'_n})$ 
to be the complement of a pair
of symplectic submanifolds of a rational surface.
Indeed, for $n\ge 3$ we always have 
${\rm vol\ } \wtilde M>{\rm vol\ } M + {\rm vol\ } B'_n$.

\begin{figure} 
\begin{center}
        \psfragscanon
        \psfrag{mnm2}{$-n-2$}
        \psfrag{m2}{$-2$}
        \psfrag{np1}{$n+1$}
        \psfrag{mnp1}{$-n+1$}
        \psfrag{nm2}{$n-2$}
        \psfrag{m4}{$-4$}
        \psfrag{m3}{$-3$}
        \psfrag{0}{$0$}
        \psfrag{1}{$1$}
        \psfrag{2}{$2$}
        \psfrag{3}{$3$}
        \psfrag{M}{$M$}
        \psfrag{CP2}{$\bC P^2$}
        \psfrag{F2}{$F_2$}
        \psfrag{F3}{$F_3$}
        \psfrag{F4}{$F_4$}
        \psfrag{Fnm2}{$F_{n-2}$}
        \psfrag{Fnm1}{$F_{n-1}$}
        \includegraphics[width=4.75in]{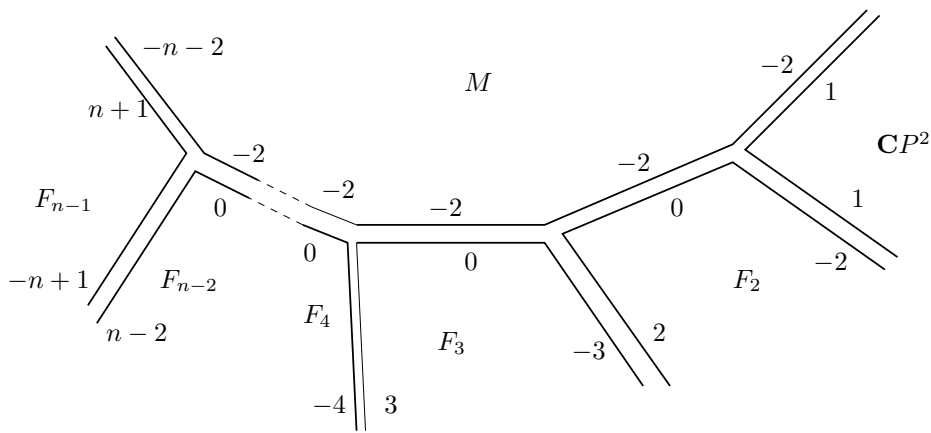}
\end{center}
\caption{Rational blowdown using the $3$-fold sum}
\label{3sum.fig}
\end{figure}

\begin{rmk} \label{3sum.rmk} \rm
This theorem can also be proved using the $3$-fold sum, an adaptation
of the symplectic sum for positively intersecting surfaces.  For details
on the $3$-fold sum, the reader should consult Symington~\cite{Sy2}.
Appealing to the $3$-fold sum, the proof can be encapsulated in
a figure.  
A blowdown is the sum of $M$, one copy of $\bC P^2$, and 
$n-2$ ruled surfaces $F_2,\ldots, F_{n-1}$.  
Figure~\ref{3sum.fig} shows how these manifolds
should be glued together.  
A $3$-fold sum is performed at each intersection point.  
Each line segment in the diagram represents a surface along which we are
gluing; the numbers labeling the edges are the
self-intersection numbers of the corresponding surfaces.
Note that the sum
of the self-intersection numbers of each pair of corresponding surfaces equals
the negative of the number of $3$-fold sums that involve the two surfaces,
as it must be to perform the sum.
\end{rmk}

\begin{rmk} \label{Lnm.rmk} \rm
Following the same method as in Example~\ref{nbhd.ex}, the
neighborhood of any union of spheres defined by a linear plumbing
and having negative definite intersection form has a symplectic
model with a symplectic structure inherited from $T^*T^2$ via
boundary reduction.
Furthermore, after removing the spheres, this model embeds into some
$(V_{n,m},\om_{n,m})$.  Indeed, if the spheres $\{S_i\}_{i=1}^s$
are indexed so that $S_i\cdot S_j=1$ if $j=i+1$ and $S_i\cdot S_j=0$
otherwise, and if they have self-intersection
numbers $-b_1,\ldots -b_s$, then the model embeds symplectically 
in $(V_{n,m},\om_{n,m})$ where $n,m$ are the relatively prime
positive integers such that
$\frac{n}{m}=[b_1,\ldots,b_s].$

It is worth noting that this embedding shows that inside any
neighborhood of such a chain of spheres one can find 
a neighborhood with $\omega$-convex 
boundary\footnote{The boundary of a symplectic manifold $(M,\om)$ is
{\bf $\om$-convex} if there is an outward-pointing vector field $X$, defined
in a neighborhood of the boundary, such that 
$\cL_X\om=\om$.  Such an $X$ is called an {\bf expanding vector field}.
See~\cite{EG} and~\cite{Et2} for discussions of symplectic convexity.}.
Indeed, the vector field
$X=p_1 \frac{\del}{\del p_1} + p_2\frac{\del}{\del p_2}$ on $T^*T^2$
induces an expanding vector field $X_{n,m}$ on $V_{n,m}$,
and hence on the complement of the spheres $\{S_i\}_{i=1}^s$
in their model neighborhood.  
Thus, if the boundary of a neighborhood is transverse to this
expanding vector field, then it is $\omega$-convex.
\end{rmk}

\begin{question} \rm
Park~\cite{Pa} has studied, in the smooth category,
a generalized rational blowdown in which
rational balls with boundary $L(n^2,nk-1)$, $(n,k)=1$, replace a
neighborhood of a chain of spheres with the same boundary.
It would be interesting to know whether or not it is possible to
perform this generalized blowdown in the symplectic category.
\end{question}

{\it Acknowledgments}
I would like to thank Ron Fintushel for encouraging me to show that the
existence of the 
symplectic rational blowdown follows from an application of the $3$-fold sum.
Thanks to Dusa McDuff and John Etnyre for helpful comments.
Also, I am appreciative of the hospitality of
the mathematics departments at the State University of New York at Stony Brook
and the University of Arizona.

\end{document}